\documentclass[11pt]{amsart}

\usepackage{amssymb,fullpage}
\usepackage{amsrefs}
\usepackage[all]{xy}

\newcommand{\R}{\mathbb{R}}
\newcommand{\C}{\mathbb{C}}
\newcommand{\Z}{\mathbb{Z}}

\swapnumbers
\newtheorem{thm}{Theorem}[section]
\newtheorem{lemma}[thm]{Lemma}
\newtheorem{theorem}[thm]{Theorem}

\newtheorem{proposition}[thm]{Proposition}
\newtheorem*{corollary*}{Corollary}

\theoremstyle{definition}
\newtheorem{definition}[thm]{Definition}
\newtheorem{remark}[thm]{Remark}
\newtheorem{example}[thm]{Example}

\newtheorem {empt}[thm]      {}

\numberwithin{equation}{section}

\begin{document}

\title[Symplectic Deligne--Mumford stacks]{Hamiltonian group actions 
on symplectic Deligne-Mumford stacks and toric orbifolds}

\author{Eugene Lerman}
\author{Anton Malkin}
\thanks{Supported in
part by NSF grants}
\address{Department of Mathematics, University of Illinois, Urbana, IL 61801}

\begin{abstract}
We develop differential and symplectic geometry of differentiable 
Deligne-Mumford stacks (orbifolds) including Hamiltonian group actions
and symplectic reduction. As an application we construct new examples of 
symplectic toric DM stacks.
\end{abstract}

\maketitle
\tableofcontents
\section{Introduction}

We have three goals in this paper. The most fundamental is to write
down in a consistent form the basics of differential and symplectic
geometry of orbifolds thought of as Deligne-Mumford (DM) stacks over the
category of smooth manifolds. This includes descriptions of the
tangent and cotangent bundles, vector fields, differential forms, Lie group
actions, and symplectic reduction.  Most if not all of these notions are
well-known on the level of being ``analogous to
manifolds".  Recall that in the original approach of Satake 
\cite{Satake1956}  an  orbifold is a topological space which  is 
locally a quotient of a vector space by a finite group action.  Smooth
functions invariant under these local group actions form the structure
sheaf.  A more recent incarnation of this idea, largely due to
Haefliger, is to think of an atlas on an orbifold as a proper etale
Lie groupoid \cite{Moerd}.  This approach makes it easy to define
local geometric structures such as vector fields, differential forms,
symplectic structures and Morse functions. However global structures
such as Lie group actions are awkward to work with in an \'etale
atlas.  One of our observations is that global structures look much
simpler in suitable non-\'{e}tale atlases. So we prefer to think of
orbifolds as a Deligne-Mumford (DM) stacks and compute in arbitrary
atlases.  The downside is that in an arbitrary groupoid atlas vector
fields and differential forms look more complicated.  We show that
there are consistent descriptions of all such geometric structures on
a DM stack.  More specifically, given a DM stack $\mathcal{X}$ there
is a presentation $X_1 \rightrightarrows X_0 \rightarrow
\mathcal{X}$ of $\mathcal{X}$ by a Lie groupoid $X_1 \rightrightarrows
X_0$ so that any geometric structure on $\mathcal{X}$ is given by a
compatible pair of the corresponding structures on $X_1$ and $X_0$.
For example a vector field (differential form, function) is a
compatible pair of vector fields (differential forms, functions) on
$X_1$ and $X_0$. Similarly given a Lie group action on $\mathcal{X}$
there is an atlas $X_0\to \mathcal{X}$ so that the action can be
described by a pair of \emph{free} actions on $X_1$ and $X_0$. Such a
presentation of a group action is useful even in the case of
manifolds, where it can be thought of as a stacky version of replacing
a $G$-manifold $M$ with $EG \times_{_G} M$.  Consequently the
quotient of $\mathcal{X}$ with respect to a $G$-action is represented
by the quotients of $X_1$ and $X_0$.  Similar statements hold for
symplectic quotients, etc.  A reader not comfortable
with the abstract stack theory can safely take these pair-based
descriptions as definitions. This is perfectly fine for applications,
since the actual calculations are always done in atlases.  However to
show that the definitions make sense one should either prove that they
are atlas-independent (i.e., Morita-invariant) or convince oneself
that there is an abstract definition in terms of the stack
$\mathcal{X}$, the approach taken in this paper.

The paper is organized as follows.  In Section~\ref{sec:1} we discuss
vector fields and forms on DM stacks and provide their description in
non-\'etale atlases.  We end the section with a definition of a
symplectic DM stack.

In Section~\ref{sec:2} we review
group actions on stacks following Romagny \cite{Romagny2005}, define
Hamiltonian actions and prove an analogue of Marsden-Weinstein-Meyer
reduction theorem for DM stacks.  

In Section~\ref{sec:3} we relate group actions on quotient stacks to
group extensions. We then describe its symplectic analogue, which may be
thought of as the stacky version of reduction in stages.  

In Section~\ref{sec:4} we take up symplectic toric DM stacks.  Recall
that symplectic toric manifolds are analogues of toric varieties in
algebraic geometry, though symplectic-algebraic correspondence is not
1-1.  Compact connected symplectic toric manifolds were classified by
Delzant \cite{Delzant1988}.  Delzant's classification was extended to
compact orbifolds by Lerman and Tolman
\cite{LermanTolman1997}. However the class of orbifolds is not as
natural as the class of DM stacks. For example it is not closed under
taking substacks.
For this reason we feel it is preferable to work with symplectic toric
DM stacks rather than orbifolds. 

In algebraic geometry the corresponding notion of a toric  DM
stack is still evolving.  To the best of our knowledge, it first
appeared  in the work of Borisov, Chen, and Smith
\cite{BorisovChenSmith2005} as a construction. 
Later Iwanari \cite{Iwanari2006} proposed the definition 
of a toric triple as an
effective DM stack with an action of an algebraic torus having a dense
orbit isomorphic to the torus.  Recently, Fantechi, Mann and Nironi
\cite{FantechiMannNironi2007} gave a new definition of a smooth toric DM stack 
as DM stack with an action of a DM torus $\mathcal{T}$ having a dense
open orbit isomorphic to $\mathcal{T}$. According to 
\cite{FantechiMannNironi2007}, a DM torus is a Picard
stack isomorphic to $T \times B \Gamma$ where $T$ is an 
algebraic torus and $B\Gamma$ is the classifying stack of 
a finite \emph{abelian} group $\Gamma$. 

We define a symplectic toric DM stack as a symplectic DM stack with an
effective Hamiltonian action of a compact torus. Then, generalizing a
construction in \cite{Delzant1988}, we produce a large class of
examples of symplectic toric DM stacks as symplectic quotients of the
form $(\mathbb{C}^N \times B\Gamma)/\!\!/_c\, A$, where $\Gamma$ is an
\emph{arbitrary} finite group and $A$ is a closed subgroup of
$\R^N/\Z^N$. From the point of view of symplectic geometry the
restriction that $\Gamma$ is abelian is unnatural.  Note additionally
that though we start with a trivial gerbe $\mathbb{C}^N \times
B\Gamma$ over $\mathbb{C}^N$, the resultant toric DM stack
$(\mathbb{C}^N
\times B\Gamma)/\!\!/_c \,A$ may be nontrivial as a gerbe over
$\mathbb{C}^N /\!\!/_c \,A$.

\section{Differential and symplectic geometry of DM stacks}\label{sec:1}

\begin{empt}\textbf{Groupoids and stacks.}\quad
All stacks in this paper are stacks over the category of smooth manifolds
with the submersion Grothendieck topology
\citelist{\cite{BehrendXu2006} \cite{Metzler2003} \cite{Heinloth2005}}.
Recall that a stack $\mathcal{X}$ is \emph{differentiable} if 
there is an atlas (representable surjective submersion) 
$X_0 \rightarrow \mathcal{X}$, where $X_0$ is a manifold.
Given an atlas one has a \emph{presentation}
$X_1 = X_0 \times_{_\mathcal{X}} X_0 \rightrightarrows X_0 
\rightarrow \mathcal{X}$
of $\mathcal{X}$ by a Lie groupoid.  We use the notation $X_1
\rightrightarrows X_0$ for a Lie groupoid with the space of objects
$X_0$ and the space of arrows $X_1$ (cf. \cite{MoerdijkMrcun2003}).
The two maps from $X_1$ to $X_0$ are the \emph{source} and
\emph{target} maps denoted by $s$ and $t$ respectively; we suppress
the rest of the structure maps of the groupoid.  Different
presentations of the same stack are Morita equivalent as groupoids
(\citelist{\cite{BehrendXu2006} \cite{Metzler2003}
\cite{Heinloth2005}}).  So one can think of differentiable stacks,
to first approximation, as Morita-equivalence classes of Lie
groupoids (this doesn't quite work when group actions enter the picture).

The main goal of this section is to describe various geometric
strictures on DM stacks. As mentioned in the introduction, the common
feature of these descriptions is that there is always a presentation
$X_1 \rightrightarrows X_0 \rightarrow \mathcal{X}$ in which the
structure on $\mathcal{X}$ is given by a compatible pair of the
corresponding structures on $X_1$ and $X_0$. A reader not comfortable
with the abstract stack theory can safely take these pair-based
descriptions as definitions. This is perfectly fine for applications,
since the actual calculations are always done in atlases.  However to
show that the definitions make sense one should either prove that they
are atlas-independent (i.e., Morita-invariant) or convince oneself
that there is an abstract definition in terms of the stack
$\mathcal{X}$.

Recall that a stack $\mathcal{X}$ is \emph{Deligne-Mumford} (DM) if
there is a presentation $U_1 \rightrightarrows U_0 \rightarrow
\mathcal{X}$ such that the groupoid $U_1 \rightrightarrows U_0$ is
\'{e}tale (i.e., $s$ and $t$ are local diffeomorphisms) and
proper (i.e., the map $s \times t :U_1 \rightarrow U_0 \times U_0$ is
proper).  Note that even if $\mathcal{X}$ is DM it is often useful and
sometimes necessary to consider non-\'{e}tale atlases
of $\mathcal{X}$. See examples below.
\end{empt}
\begin{remark}
The classical Satake definition of an orbifold can be
reformulated as that of a DM stack with an \'etale atlas
$U_1\rightrightarrows U_0$ such that the stabilizers in $U_1$ of points
in an open dense subset of $U_0$ are trivial.  This condition is not
preserved under taking substacks. So we consider a more
general notion of a DM stack.  In such a stack every point could have
a non-trivial finite stabilizer group (also called \emph{inertia}
group). An example is the classifying stack $B\Gamma$ of a finite
group $\Gamma$, which can be presented by the groupoid $\Gamma
\rightrightarrows \mathrm{pt}$.
\end{remark}

\begin{empt} \textbf{Lie algebroids and 
non-\'{e}tale presentations of DM stacks.}\quad
Recall (cf. \cite{MoerdijkMrcun2003}) that the Lie algebroid of a Lie
groupoid $X_1 \rightrightarrows X_0$ is a vector bundle map called
the anchor $A \xrightarrow{a} TX_0$, where $A = \ker ds |_{_{X_0}}$,
and $a=dt$.

The following theorem describes DM stacks in terms of their arbitrary
(not necessarily \'{e}tale) atlases. Such atlases
will be crucial later when we study Lie group actions 
on DM stacks.
\end{empt}
\begin{thm}\label{DMDefinitionsTheorem}
Let $\mathcal{X}$ be a stack over the category of
smooth manifolds. The following conditions are equivalent:
\begin{enumerate}
\item
$\mathcal{X}$ is DM, i.e., has a proper \'{etale} presentation
$U_1 \rightrightarrows U_0 \rightarrow \mathcal{X}$;
\item
$\mathcal{X}$ has a proper presentation 
$X_1 \rightrightarrows X_0 \rightarrow \mathcal{X}$
such that the anchor map of the corresponding Lie algebroid
is injective;
\item
$\mathcal{X}$ has a  proper presentation 
$X_1 \rightrightarrows X_0 \rightarrow \mathcal{X}$
such that the subbundles 
$\ker ds$ and $\ker dt$ of $TX_1$ are transverse;
\item
$\mathcal{X}$ is proper differentiable (has a proper atlas) and
for any presentation 
$X_1 \rightrightarrows X_0 \rightarrow \mathcal{X}$
the anchor map of the corresponding Lie algebroid is injective;
\item
$\mathcal{X}$ is proper differentiable and for any presentation 
$X_1 \rightrightarrows X_0 \rightarrow \mathcal{X}$
the subbundles 
$\ker ds$ and $\ker dt$ of $TX_1$ are transverse;
\item
$\mathcal{X}$ is proper differentiable (has a proper atlas) and
the inertia group of any point of $\mathcal{X}$ is finite.
\end{enumerate}
\end{thm}
\begin{proof} This is essentially Theorem 1 of \cite{CrainicMoerdijk2001}.
We sketch the main ideas of the proof. To go from an \'{e}tale
presentation $U_1 \rightrightarrows U_0 \rightarrow \mathcal{X}$ to an
arbitrary presentation $X_1 \rightrightarrows X_0 \rightarrow
\mathcal{X}$ consider the pullback diagram
\begin{equation}\nonumber
\xymatrix{
U \times_{\scriptscriptstyle \mathcal{X}} X_1 
\ar@<0.6ex>[r] \ar@<-0.2ex>[r] \ar[d] &
U \times_{\scriptscriptstyle \mathcal{X}} X_0 \ar[r] \ar[d] &
U \ar[d] \\
X_1 \ar@<0.6ex>[r] \ar@<-0.2ex>[r] &
X_0 \ar[r] &
\mathcal{X}
}
\end{equation}
Note that the top row is a presentation of the manifold
$U$ and hence the anchor map of the corresponding algebroid
is injective. Now, since the vertical maps are \'{e}tale,
the same holds for the bottom row algebroid.

In the opposite direction, given a presentation
$X_1 \rightrightarrows X_0 \rightarrow \mathcal{X}$
such that inertia groups are all finite or
equivalently the corresponding Lie algebroid 
has injective anchor map, the action of $X_1$ defines
a foliation of $X_0$.  One  constructs an \'{e}tale
atlas of $\mathcal{X}$ using transverse slices to this
foliation. We refer the reader to  \cite{CrainicMoerdijk2001}
for  details. 
\end{proof}

\begin{empt}\textbf{Differential forms and vector fields.}\quad
A smooth manifold $X$ comes equipped with the tangent $TX$ and the
cotangent $T^*X$ bundles, the sheaf of vector fields
$\mathrm{Vect}_X$, and the de Rham complex of sheaves of differential
forms $\Omega_X^{\bullet}$.  Moreover $\mathrm{Vect}_X$ is the sheaf
$\mathcal{C}^{\infty}_{TX}$ of smooth sections of $TX \rightarrow X$,
and $\Omega_X^{\bullet}$ is the sheaf
$\mathcal{C}^{\infty}_{\Lambda^{\bullet}T^*X}$ of smooth sections of
$\Lambda^{\bullet} T^*X \rightarrow X$. The only part of this story
compatible with pullbacks and hence defined for stacks is the de Rham
complex. Namely, given an arbitrary stack $\mathcal{X}$ the de Rham
complex of sheaves on $\mathcal{X}$ is defined as follows (see, for
example, \cite{BehrendXu2006}): for an object $\upsilon \in
\mathcal{X}$ over a manifold $U$ (in other words, a map $\upsilon: U
\rightarrow \mathcal{X}$) one has $\Omega_{\mathcal{X}}^{\bullet}
(\upsilon) = \Omega_U^{\bullet} (U)$, the de Rham complex on $U$.  A
differential form of degree $k$ on $\mathcal{X}$ is a global section
of the sheaf $\Omega^k_{\mathcal{X}}$, i.e., a homomorphism from the
trivial sheaf on $\mathcal{X}$ to $\Omega^k_{\mathcal{X}}$.  However
the sheaf $\Omega_{\mathcal{X}}^k$ is not the sheaf of sections of a
vector bundle on $\mathcal{X}$ even if $\mathcal{X}$ is a DM stack.
Moreover, though one can define tangent stack $T\mathcal{X}$ for an
arbitrary stack $\mathcal{X}$ \cite{Hepworth2008}, the projection
$\pi: T\mathcal{X} \rightarrow \mathcal{X}$ is not a vector
bundle. Rather it is a 2-vector bundle. Hence sections of $\pi$ form
not a sheaf of sets but a sheaf of groupoids. We refer the reader to
\cite{BehrendXu2006} and \cite{Hepworth2008} for discussions on 
sections of a map of stacks. 

The situation is much better in the case of DM stacks, which is the
reason na\"{i}ve definitions used in orbifold theory work well for
many purposes. However since we need to use arbitrary and not only
\'{e}tale atlases for DM stacks,  we explain the
concepts associated with tangent and cotangent bundles and their
sections in some detail.  A reader familiar with foliations or locally
free group actions should recognize many constructions, such as
presentations of transverse tangent bundles via Lie algebroids.

The crucial properties of an \'{e}tale map (local diffeomorphism) 
$f : M \rightarrow N$ is that
one can pull-back vector fields along $f$ and that the pull-back 
of the (co)tangent bundle is the (co)tangent bundle 
(i.e., $f^*(TN) = TM$, $f^* (T^*N) = T^*M$).
Hence the following definitions make sense for a DM stack $\mathcal{X}$.
Consider an \'{e}tale presentation 
$U_1 \rightrightarrows U_0 \rightarrow \mathcal{X}$.
Then we have $s^* (TU_0) = t^* (TU_0) = TU_1$ and hence the bundle 
$TU_0 \rightarrow U_0$ descends to a vector bundle 
$T \mathcal{X} \rightarrow \mathcal{X}$ called the tangent bundle of
$\mathcal{X}$. It is easy to see
that the bundle $T \mathcal{X} \rightarrow \mathcal{X}$ does not depend 
(up to equivalence) on the choice of the \'{e}tale atlas. One defines
the cotangent bundle $T^{\ast} \mathcal{X} \rightarrow \mathcal{X}$
in a similar way. These are the usual definitions in the orbifold
theory phrased in a fancy way. 

Let us now consider an arbitrary presentation 
$X_1 \rightrightarrows X_0 \xrightarrow{\xi} \mathcal{X}$ of our
DM stack $\mathcal{X}$. We would like to describe pullbacks 
$\xi^* (T\mathcal{X})$ and $\xi^* (T^*\mathcal{X})$ of the 
tangent $T\mathcal{X}$ and the cotangent bundle $T^* \mathcal{X}$. 
Let $A \xrightarrow{a} TX_0$ be the Lie algebroid
of the Lie groupoid $X_1 \rightrightarrows X_0$, that is 
$A = \ker ds |_{_{X_0}}$, $a=dt$.
Since $\mathcal{X}$ is DM, Theorem 
\ref{DMDefinitionsTheorem} implies that 
$\ker ds$ and $\ker dt$ are transverse as subbundles of $TX_1$ and 
the anchor map $a$ is injective. We identify $A$ with its image in $TX_0$ 
in what follows.
\end{empt}
\begin{proposition}
\label{ArtinBundleProposition}
Let $X_1 \rightrightarrows X_0 \xrightarrow{\xi} \mathcal{X}$
be a presentation of a DM stack $\mathcal{X}$ with the associated
algebroid $A \hookrightarrow TX_0$. Then
\begin{enumerate}
\item\label{ArtinBundle1}
$\xi^* (T\mathcal{X}) = TX_0/A$ as a vector bundle on $X_0$;
\item
$s^* ( \xi^* (T\mathcal{X})) = 
t^* ( \xi^* (T\mathcal{X})) = TX_1/(\ker ds + \ker dt)$ 
as a vector bundle on $X_1$;
\item
$\xi^* (T^*\mathcal{X}) = A^{\perp}$, where 
$A^\perp \subset T^*X_0 = (TX_0)^*$ is the annihilator of $A$;
\item\label{ArtinBundle4}
$s^*(\xi^* (T^*\mathcal{X})) = t^*(\xi^* (T^*\mathcal{X})) = 
(\ker ds )^{\perp} \cap (\ker dt )^\perp \subset T^*X_1$.
\end{enumerate}
\end{proposition}
\begin{proof}
Let $U \rightarrow \mathcal{X}$ be an \'{e}tale atlas of $\mathcal{X}$.
Then in the diagram 
\begin{equation}\nonumber
\xymatrix{
U \times_{\scriptscriptstyle \mathcal{X}} X_1 
\ar@<0.6ex>[r] \ar@<-0.2ex>[r] \ar[d] &
U \times_{\scriptscriptstyle \mathcal{X}} X_0 \ar[r] \ar[d] &
U \ar[d] \\
X_1 \ar@<0.6ex>[r] \ar@<-0.2ex>[r] &
X_0 \ar[r] &
\mathcal{X}
}
\end{equation}
vertical maps are \'{e}tale. Hence it is enough to prove 
\eqref{ArtinBundle1}--\eqref{ArtinBundle4}
for the first row (a presentation of the manifold $U$), which a standard
exercise in differential geometry of fibrations.
\end{proof}

%\begin{empt}\textbf{Differential forms and vector fields on DM stacks.}\quad

\begin{definition}We define sheaves of vector fields and differential forms 
on a DM stack $\mathcal{X}$ as sheaves
of smooth sections $\mathcal{C}^{\infty}_{T\mathcal{X}}$ and 
$\mathcal{C}^{\infty}_{\Lambda^{\bullet}T^*\mathcal{X}}$
of the tangent and the exterior powers of the cotangent bundle
respectively.
\end{definition}
Sections of these vector bundles do form sheaves of 
vector spaces (as opposed to sheaves of groupoids) 
as follows from the following explicit description (cf. \cite{Hepworth2008}).
The space of sections
$\mathcal{C}^{\infty}_{T\mathcal{X}} (\upsilon )$
on an \'{e}tale map $U \xrightarrow{\upsilon} \mathcal{X}$, i.e.,
on an object $\upsilon$ of $\mathcal{X}$ over $U$, is just
$\mathcal{C}^{\infty}_{TU} (U)$ -- the space of vector fields on $U$. 
For a pullback $f$ of \'{e}tale maps
\begin{equation}\nonumber
\xymatrix{
U \ar[r] \ar[d]_{f} & \mathcal{X} \\
V \ar[ur]
}
\end{equation}
the pullback $f^*: \mathcal{C}^{\infty}_{T\mathcal{X}} (V) \rightarrow \mathcal{C}^{\infty}_{T\mathcal{X}} (U)$
of sections is the pullback of vector fields under the \'{e}tale map $f$. 
Similarly
$\mathcal{C}^{\infty}_{\Lambda^{\bullet} T^*\mathcal{X}} (\upsilon )= \mathcal{C}^{\infty}_{\Lambda^{\bullet} T^*\mathcal{X}} (U)=
\mathcal{C}^{\infty}_{\Lambda^{\bullet} T^*U} (U)=
\Omega^{\bullet}_{U} (U) = \Omega^{\bullet}_{\mathcal{X}} (U)$
for an \'{e}tale map $U \xrightarrow{\upsilon} \mathcal{X}$.
\begin{definition}
We define vector fields and differential forms on a stack
$\mathcal{X}$ as global sections of the corresponding sheaves. 
\end{definition}
Recall that the vector space of global sections of a sheaf $F$ of
vector spaces on $\mathcal{X}$ is the vector space of homomorphisms
from the trivial sheaf $1_{_\mathcal{X}}$ to $F$.  Thus, given a
presentation $X_1 \rightrightarrows X_0 \xrightarrow{\xi}
\mathcal{X}$, the space of global sections of $F$ is the equalizer of
the two pull-back  maps $F(X_0) \rightrightarrows
F(X_1)$.

The following proposition describes vector fields and
differential forms in an atlas.  
Recall that given a surjective submersion 
$f: Y \rightarrow X$ of manifolds, a vector field
$v_{_Y} \in \mathrm{Vect} (Y)$ and a vector field 
$v_{_X} \in \mathrm{Vect} (X)$, 
one says that $v_{_Y}$ is \emph{$f$-related} to $v_{_X}$ if 
\[
df (v_{_Y}) = v_{_X} \circ f.
\] 
Note that this is a relation, not a map $\mathrm{Vect} (Y) \rightarrow
\mathrm{Vect} (X)$.  Given $v_{_X}$, $v_{_Y}$ is determined up to
a section of the bundle $\ker df \subset TY$.
%\end{empt}
\begin{proposition}\label{prop:2.8}
Let $X_1 \rightrightarrows X_0 \xrightarrow{\xi} \mathcal{X}$
be a presentation of a DM stack $\mathcal{X}$ with the associated
algebroid $A \hookrightarrow TX_0$. Then
\begin{enumerate}
\item\label{ArtinGlobalVectors}
The Lie algebra $\mathrm{Vect} (\mathcal{X}) :=
\mathcal{C}^{\infty}_{T\mathcal{X}} (\mathcal{X})$ of vector fields on
$\mathcal{X}$, i.e., of global sections of $T\mathcal{X}$, is isomorphic
to the Lie algebra $\mathcal{C}^{\infty}_{TX_0/A} (X_0)^{X_1}$ of
$X_1$-invariant sections of the bundle $TX_0/A$.  Explicitly, vector
fields on $\mathcal{X}$ are equivalence classes of pairs consisting of
a vector field $v_0$ on $X_0$ and a vector field $v_1$ on $X_1$, which
are both $s$- and $t$-related:
\[
\mathrm{Vect} (\mathcal{X}) \cong
\frac{ 
\bigl\{(v_1, v_0) \in \mathrm{Vect} (X_1) 
 \times \mathrm{Vect} (X_0) \mid  
ds (v_1) = v_0 \circ s , \ dt (v_1) = v_0 \circ t \bigr\}
}{
\bigl\{(v_1, v_0) \mid
ds (v_1) = v_0 \circ s , \ dt (v_1) = v_0 \circ t ,
\ v_1 \in ( \ker ds + \ker dt ) \bigr\} 
}
\]
\item\label{ArtinGlobalForms}
The de Rham complex $\Omega^{\bullet} ( \mathcal {X} ) :=
\mathcal{C}^{\infty}_{\Lambda^{\bullet} T^* \mathcal{X}} (\mathcal{X})$ 
of differential forms on $\mathcal{X}$, i.e., of global sections of $\Lambda^{\bullet} T^* \mathcal{X}$, is isomorphic to the
complex $\mathcal{C}^{\infty}_{\Lambda^* A^{\perp}} (X_0)^{X_1}$
of $X_1$-invariant forms on $X_0$:
\[
\Omega^{\bullet} (\mathcal{X}) \cong
\{ \tau \in \Omega^{\bullet} (X_0) \ | 
\ s^* \tau = t^* \tau \} \quad ,
\]
which can be also expressed as a set of pairs:
\[
\Omega^{\bullet} (\mathcal{X}) \cong
\bigl\{ ( \sigma_1,\sigma_0)  \in 
\Omega^{\bullet} (X_1) \times \Omega^{\bullet} (X_0)
\bigr)
\mid
 s^* \sigma_0 = \sigma_1 , \  t^* \sigma_0 = \sigma_1 \bigr\} \quad .
\]
\item
The contraction of vector fields and forms on $\mathcal{X}$ is induced 
by the contraction of vector fields and forms on $X_0$ and $X_1$:
\begin{equation}\nonumber
\iota_{( v_1 , v_0 )} \, ( \sigma_1 , \sigma_0 ) = 
( \iota_{v_1}  \sigma_1 , \iota_{v_0} \sigma_0 )
\end{equation}
\end{enumerate}
\end{proposition}
\begin{proof} By definition of the sheaf of sections of 
a vector bundle, pullbacks of sections are sections of pullback.
Now the proposition follows from the explicit description 
\ref{ArtinBundleProposition} of the pull-backs
of the tangent and the cotangent bundles to $X_0$. 
\end{proof}

\begin{remark}One can consider  Proposition~\ref{prop:2.8} above as a 
definition of vector fields and differential forms on a DM stack
$\mathcal{X}$. Abstract definitions in terms of global sections of
vector bundles on $\mathcal{X}$ ensure atlas-independence.  One can
also check the atlas-independence directly, without any reference to
stacks.
\end{remark}
\begin{remark} Let $X_0 \rightarrow \mathcal{X}$ be an arbitrary
(as opposed to \'{e}tale) submersion. Then sections on $X_0$ of the sheaf $\mathcal{C}^{\infty}_{\Lambda^{\bullet} T^*\mathcal{X}}$
are differential forms on $X_0$ vanishing on the corresponding 
Lie algebroid $A \hookrightarrow TX_0$.
Now recall that sections of the abstract de Rham complex $\Omega^{\bullet}_{\mathcal{X}}$ (defined for
arbitrary stacks, cf. \cite{BehrendXu2006}) are arbitrary forms on $X_0$.
Hence $\mathcal{C}^{\infty}_{\Lambda^{\bullet} T^*\mathcal{X}}$
and $\Omega^{\bullet}_{\mathcal{X}}$ \emph{are not isomorphic as sheaves.} 
In particular $\Omega^{k}_{\mathcal{X}}$ is not a sheaf of sections 
of a vector bundle on $\mathcal{X}$. 
However the spaces of global sections (differential forms on $\mathcal{X}$)
and hypercohomolgy groups (de Rham cohomology of $\mathcal{X})$
of these two sheaves are isomorphic.
For example the space of global sections of either 
$\mathcal{C}^{\infty}_{\Lambda^{\bullet} T^*\mathcal{X}}$ or 
$\Omega^{\bullet}_{\mathcal{X}}$ is isomorphic 
to the space of differential forms $\tau$ on $X_0$ satisfying
$s^* \tau = t^* \tau$ on $X_1$.
\end{remark}
\begin{empt}\textbf{Symplectic DM stacks.}\quad
Given the above abstract definitions and explicit descriptions
of vector fields and differential forms on a DM stack,
it is now straightforward to define symplectic forms on DM stacks.

A 2-form $\omega \in \Omega^2 (\mathcal{X})$ on a DM stack $\mathcal{X}$ is
\emph{non-degenerate} if the contraction with $\omega$ induces an
isomorphism $\mathrm{Vect} ( \mathcal{X} ) \rightarrow \Omega^1
(\mathcal{X})$.  In a presentation $X_1 \rightrightarrows X_0
\rightarrow \mathcal{X}$, the 2-form $\omega$ is given (represented)
by a pair $(\omega_1, \omega_0)$, and it is non-degenerate iff $\ker
\omega_0 = A \subset TX_0$, the Lie algebroid of $X_1
\rightrightarrows X_0$ or, equivalently, iff $\ker \omega_1 = \ker ds +
\ker dt \subset TX_1$

A \emph{symplectic form} $\omega$ on a DM stack $\mathcal{X}$ is a a
non-degenerate closed 2-form $\omega \in \Omega^2 (\mathcal{X})$.
A \emph{symplectic DM stack} is a pair $(\mathcal{X} , \omega )$,
where $\mathcal{X}$ is a DM stack and $\omega$ is a symplectic form on
$\mathcal{X}$.
\end{empt}

\section{Hamiltonian group actions on symplectic DM stacks}\label{sec:2}

\begin{empt}\textbf{Group actions on stacks} (following \cite{Romagny2005}).
Group actions on stacks are more complicated than actions on manifolds
because stacks are categories and the collection of all stacks forms a
2-category.  Thus when a group acts on a stack, group elements act as
functors but the composition of functors representing two group
elements can differ from the functor representing their product by a
natural transformation. In the case of Lie group $G$ action on a
differentiable stack $\mathcal{X}$ we would like the action to be
``smooth'' so we represent it as a map $a: G \times \mathcal{X}
\rightarrow \mathcal{X}$ instead of an action homomorphism from $G
\rightarrow Aut(\mathcal{X})$.  The fact that $a$ is an action is
encoded in the 2-commutativity of the diagrams
\[
\xymatrix{
G \times G \times \mathcal{X} 
\ar[r]^-{m \times \mathrm{id}_{_{\mathcal{X}}}}  
\ar[d]_{\mathrm{id}_{_G} \times a} &
G \times \mathcal{X} \ar[d]^{a} \\
G \times \mathcal{X} \ar[r]_-{a} \ar@{=>}[ur]^{\alpha} &
\mathcal{X}
} \qquad \qquad \qquad
\xymatrix{
G \times \mathcal{X}  \ar[r]^-{a} 
\ar@{=>}<7mm,-7mm>^-{\epsilon} & 
\mathcal{X} \\
\mathcal{X} 
\ar[u]^{e_{_G} \times \mathrm{id}_{_{\mathcal{X}}}}
\ar[ur]_{\mathrm{id}_{_{\mathcal{X}}}}
}
\]
where $m$ is the multiplication map in $G$, and
$e_{_G}$ is the identity inclusion. The natural transformations
$\alpha$ and $\epsilon$, which are part of the data defining the action,  
should satisfy further compatibility conditions \cite{Romagny2005}.
The whole definition may be  thought of as describing a stack over $BG$.
\end{empt}
Given an action $a: G \times \mathcal{X} \rightarrow \mathcal{X}$ 
and an atlas 
$X_0 \rightarrow \mathcal{X}$ we consider the composition
\begin{equation}\nonumber
G \times X_0 \rightarrow 
G \times \mathcal{X} \xrightarrow{a} 
\mathcal{X}
\end{equation}
to obtain a new atlas 
$\widetilde{X_0} = G \times X_0$ of $\mathcal{X}$.
Then it follows from the definition of the action of $G$ 
on $\mathcal{X}$ that left $G$-translations on 
$G \times X$ induce free $G$-actions on
both $\widetilde{X_0}$ and 
$\widetilde{X_1} = 
\widetilde{X_0} \times_{_\mathcal{X}} \widetilde{X_0}$
and these actions commute with the structure maps of
the presentation $\widetilde{X_1} \rightrightarrows \widetilde{X_0}$ 
of $\mathcal{X}$. More precisely, $G$-stacks 
$[X_0/X_1]$ and $[\widetilde{X_0}/\widetilde{X_1}]$ 
are isomorhic via an equivalence of fibered categories involving
the natural transformations $\alpha$ and $\epsilon$
from the definition of the $G$-action on $\mathcal{X}$.
We refer the reader to 
\cite{Romagny2005} for details of this construction
(called \emph{strictification} in the \emph{loc.\ cit.})
and summarize this discussion in the following proposition.

\begin{proposition} Suppose a Lie group $G$ acts on a 
differentiable stack $\mathcal{X}$. Then $\mathcal{X}$ has 
a presentation, called a \emph{$G$-presentation},
$X_1 \rightrightarrows X_0 \rightarrow \mathcal{X}$
in which the $G$-action is given by free $G$-actions on
$X_1$ and $X_0$ compatible with structure maps of 
the groupoid $X_1 \rightrightarrows X_0$.
\end{proposition}
\begin{proof} This is essentially Proposition 1.5 of
\cite{Romagny2005}. The freeness of the actions is not stated
there but follows from the proof. 
\end{proof}

One can define a $G$-stack as a stack represented by a Lie groupoid
with free $G$-actions on the set of objects and arrows. This
is a natural definition if one thinks of $G$-stacks as stacks over
$BG$, that is,  over the site of principal  $G$-bundles.

\begin{empt}\textbf{Quotient of a stack.} \quad
Let $G$ be a compact Lie group. Given a $G$-action on a differentiable
stack $\mathcal{X}$ with a $G$-presentation $X_1 \rightrightarrows X_0
\rightarrow \mathcal{X}$, there is a differentiable quotient stack
$\mathcal{X}/G$ defined by a universal property with respect to maps
to manifolds with the trivial $G$-action \cite{Romagny2005}. A
presentation of $\mathcal{X}/G$ is given by the Lie groupoid $X_1/G
\rightrightarrows X_0/G \rightarrow \mathcal{X}/G$. Note that $X_1/G$
and $X_0/G$ are manifolds since $G$-actions on $X_0$ and $X_1$ are
free and proper.

This is the only place in the paper we use properness/compactness
conditions. One can consider non-proper DM-stacks and arbitrary
Lie group actions with the quotient stack 
represented by the semidirect product groupoid 
$G \times X_1 \rightrightarrows X_0$.
We would like to stick to our philosophy that every structure/procedure
(in particular, a $G$-action and quotient) 
should be represented by a pair of
the corresponding structures/procedures on objects and arrows 
of an appropriate presentation. So we restrict ourselves to
proper actions.
\end{empt}
\begin{example} 
Let a compact Lie group $G$ act on the classifying stack $BH$ of a
compact Lie group $H$. The stack $BH$ has a presentation $H
\rightrightarrows \mathrm{pt} \rightarrow BH$. Hence a
$G$-presentation of $BH$ is given by $X_0 = G \times \mathrm{pt} =G$,
$X_1 = G \times_{_{BH}} G = K \times G$, where $K$ is a principal
$H$-bundle over $G$ and also a group. Therefore $G$-actions on $BH$
correspond to Lie group extensions $1 \rightarrow H \rightarrow K
\rightarrow G \rightarrow 1$.  The $G$-presentation of $BH$
corresponding to such an extension is given by the action groupoid $K
\times G \rightrightarrows G \rightarrow BH$ for the right action of
$K$ on $G$, and $G$-action on $BH$ in this presentation comes from the
left action of $G$ on $G$. The quotient stack $BH/G$ is equivalent to
$BK$.  See  \ref{MTimesHSection} for a generalization of this example 
and of the example below.
\end{example}
\begin{example}
Suppose a Lie group $G$ acts on a manifold $M$. Then we
have a $G$-atlas $X_0 = M \times G$ of $M$, 
with the $G$-presentation given by the action groupoid
$G \times (M \times G) \rightrightarrows M \times G \rightarrow M$
corresponding to the following action of $G$ on $M \times G$:
$g \cdot (m,g') = (g \cdot m, g'g^{-1})$. 
The action of $G$ on $M$ in this presentation
is given by left translations on $G$: $g \cdot (m,g') = (m, gg')$.
Hence, even in the case of manifolds, stacky point of view
has advantages: one can replace arbitrary action on a manifold by a 
free action on a groupoid (i.e., ``resolve" of the original action).
The quotient stack is $M/G$ with the presentation
$G \times M \rightrightarrows M \rightarrow M/G$.
\end{example}
\begin{empt}\textbf{Infinitesimal actions.}\quad
Given an action $a : G \times \mathcal{X} \rightarrow \mathcal{X}$
of a Lie group $G$ on a DM stack $\mathcal{X}$, we obtain the 
derived map (infinitesimal action) 
\begin{equation}\nonumber
d a: (\mathfrak{g}, 0) 
\hookrightarrow 
\mathrm{Vect} (G \times \mathcal{X})
\rightarrow 
\mathrm{Vect} (\mathcal{X}) 
\quad ,
\end{equation}
where we think of the Lie algebra $\mathfrak{g}$ as the space of
right-invariant vector fields on $G$. Moreover, though we don't use it
in this paper, infinitesimal actions on DM stacks have all the usual
properties of infinitesimal actions on manifolds. For example, they
are homomorphisms of Lie algebras. The proofs are identically the same
as in the manifold case since in the case of DM stacks the natural
transformations involved in the definition of the action act trivially
on tangent spaces and their maps. Roughly speaking, it is hard to
describe a Lie group action on an orbifold, but a Lie algebra action
is just given by  vector fields.

Now suppose $X_1 \rightrightarrows X_0 \rightarrow \mathcal{X}$
is a $G$-presentation of $\mathcal{X}$.
Differentiating free actions of $G$ on $X_1$ and $X_0$
we get a pair of vector fields 
$(v_1 (\varepsilon ), v_0 (\varepsilon ))\in \mathrm{Vect} ( X_1 ) \times 
\mathrm{Vect} ( X_0 ))$ for every 
$\varepsilon \in \mathfrak{g}$. Moreover, since the structure maps of
$X_1 \rightrightarrows X_0$ commute with the $G$-actions on $X_1$ and
$X_0$, we have $ds (v_1 (\varepsilon )) = v_0 (\varepsilon ) \circ s$
and $dt (v_1 (\varepsilon )) = v_0 (\varepsilon ) \circ t$.  Hence
(\emph{cf.} \ref{ArtinGlobalVectors}) the pair $(v_1 , v_0)$ defines a
vector field on the stack $\mathcal{X}$.  This vector field is $d a
(\epsilon)$ in the presentation $X_1 \rightrightarrows X_0$. One can
consider this description as the definition of $d a
(\varepsilon)$. The abstract definition  above it ensures
atlas-independence.
\end{empt}
\begin{empt}\textbf{Locally free actions.}\quad
Similar to the manifold case we say that an action
$a : G \times \mathcal{X} \rightarrow \mathcal{X}$
of a Lie group $G$ on a DM stack $\mathcal{X}$ is \emph{locally free} 
if the corresponding action of the Lie algebra is free, 
i.e., $d a (\varepsilon )$ is a nonvanishing vector
field for every $0 \neq \varepsilon \in \mathfrak{g}$. 
One should not confuse this condition with
the freeness of $G$ action in an atlas. For example,
in a $G$-presentation 
$X_1 \rightrightarrows X_0 \rightarrow \mathcal{X}$
the group $G$ acts freely on both $X_1$ and $X_0$, and
such a presentation exists for arbitrary $G$-action on $\mathcal{X}$.
\end{empt}
\begin{lemma}\label{InfinitesimalAction}
Let $a : G \times \mathcal{X} \rightarrow \mathcal{X}$
be an action of a compact Lie group $G$ on a DM stack $\mathcal{X}$,
with a $G$-presentation $X_1 \rightrightarrows X_0 \rightarrow \mathcal{X}$.
In particular, we have free $G$-actions on $X_1$ and $X_0$
and the quotient stack $\mathcal{X}/G$ is represented by
$X_1/G \rightrightarrows X_0/G \rightarrow \mathcal{X}/G$.
Let $B_1 \subset T X_1$ and $B_0 \subset TX_0$ be the subbundles 
spanned by infinitesimal vector fields generating the (free) actions of the
Lie algebra $\mathfrak{g}$ on $X_1$ and $X_0$ respectively,
and $A \hookrightarrow TX_0$ the Lie algebroid of 
$X_1 \rightrightarrows X_0$.
Then the following conditions are equivalent:
\begin{enumerate}
\item\label{LocallyFreeByDef}
The $G$-action on $\mathcal{X}$ is locally free; 
\item\label{LocallyFreeTransverse1}
The subbundles $B_1$, $\ker ds$ and $\ker dt$, of $TX_1$ are transverse;
\item\label{LocallyFreeTransverse0}
The subbundles $B_0$ and $A$ of $TX_0$ are transverse;
\item\label{LocallyFreeDM}
The quotient stack $\mathcal{X}/G$ is DM.
\end{enumerate}
\end{lemma}
\begin{proof}
The explicit description \eqref{ArtinGlobalVectors} of vector fields
on a differentiable stack implies that conditions 
\eqref{LocallyFreeByDef}, \eqref{LocallyFreeTransverse1}, 
and \eqref{LocallyFreeTransverse0}
are equivalent. The source $[s]$ and the target $[t]$ maps
of the quotient groupoid $X_1/G \rightrightarrows X_0/G$ are induced 
by the source and target maps of the original groupoid 
$X_1 \rightrightarrows X_0$. Hence $\ker[ds]$ and $\ker[dt]$ are
transverse (equivalently,  $\mathcal{X}/G$ is DM, cf. 
\eqref{DMDefinitionsTheorem}) iff 
\eqref{LocallyFreeTransverse1} is satisfied.
\end{proof}

\begin{empt}\textbf{DM stacks given by equations.}\quad
Let $\mathcal{X}$ be a DM stack, $f: \mathcal{X} \rightarrow
\mathbb{R}^n$ a function, and
$\xi \in \mathbb{R}^n$.  Consider the substack $f^{-1} ( \xi )$. In
general this is not a differentiable stack.  Let $X_1
\rightrightarrows X_0 \rightarrow \mathcal{X}$ be a presentation of
$\mathcal{X}$. Then $f = (f_1, f_0)$, where $f_0 : X_0 \rightarrow
\mathbb{R}^n$ and $f_1 = f_0 \circ s = f_0 \circ t : X_1 \rightarrow
\mathbb{R}^n$.  We say that $\xi$ is a \emph{regular value} of $f$ if
it is a regular value of $f_0$, that is, if the differential
$(df_0)_x$ is surjective at any point $x \in {f_0}^{-1}(\xi)$.  Note
that the surjectivity condition on $df_0$ is preserved under
precomposition with submersions. Hence the regularity condition on $f$
does not depend on the choice of a presentation of $\mathcal{X}$.  And
if $\xi$ is a regular value of $f_0$ then it is a regular value of
$f_1$.  

Assume $\xi$ is a regular value of $f: \mathcal{X} \rightarrow
\mathbb{R}^n$.  Then $f_1^{-1} ( \xi ) \rightrightarrows f_0^{-1} (
\xi )$ is a Lie subgroupoid of $X_1 \rightrightarrows X_0$
representing $f^{-1} ( \xi )$. Hence $f^{-1} ( \xi )$ is a
differentiable stack.  If we assume that $\mathcal{X}$ is DM then
$f^{-1} ( \xi )$ is also DM (for example, by the above argument in an
\'{e}tale presentation of $\mathcal{X}$).  We  record these observations 
as a lemma.
\end{empt}
\begin{lemma}\label{RegularValue}
Let $\mathcal{X}$ be a differentiable stack, $f: \mathcal{X}
\rightarrow \mathbb{R}^n$ be a function, and $\xi$ a regular value of
$f$. Then $f^{-1} ( \xi )$ is a differentiable stack. If $\mathcal{X}$
is DM then so is $ f^{-1} ( \xi )$.
\end{lemma}
\begin{empt}\textbf{Hamiltonian actions.}\quad
Let $\mathcal{X}$ be a DM stack, $\omega \in \Omega^2 (\mathcal{X})$
a closed 2-form, and $a : G \times \mathcal{X} \rightarrow \mathcal{X}$
an action of a Lie group $G$ on $\mathcal{X}$. If $\omega$ is 
non-degenerate then $(\mathcal{X} , \omega)$ is a symplectic DM stack, 
but this condition is not important for the following definition.

We say that an action $a : G \times \mathcal{X} \rightarrow \mathcal{X}$
is \emph{Hamiltonian} if there exists an equivariant  map
$\mu: \mathcal{X} \to \mathfrak{g}^*$ such that
\begin{equation}\label{MomentStack}
\iota_{_{d a (\varepsilon )}} \omega = 
d \langle \varepsilon, \mu \rangle
\end{equation}
for any $\varepsilon \in \mathfrak{g}$.
We refer to $\mu$ as a \emph{moment map}.

In a $G$-presentation 
$X_1 \rightrightarrows X_0 \rightarrow \mathcal{X}$
there are two points of view on differential forms. 
If one thinks of $\omega$ and $\mu$ as forms on $X_0$
satisfying pull-back conditions on $X_1$ then 
the equation \eqref{MomentStack} reads
\begin{equation}\label{MomentAtlas}
\iota_{_{u_\varepsilon}} \omega = 
d \langle \varepsilon, \mu \rangle
\end{equation}
where $u_\varepsilon$ is the $\mathrm{Vect} (X_0)$-component
of the vector field $d \gamma (\varepsilon)$. Essentially,
in a $G$-presentation we have a Hamiltonian action 
on $X_0$. Note that, while 
$u_\varepsilon$ is defined only up to addition of
sections of the algebroid $A$, 
this ambiguity does not matter in the equation
\eqref{MomentAtlas} because $\omega$ vanishes on $A$.

If one thinks of differential forms, vector fields, actions, etc., as
compatible pairs of the corresponding objects on $X_0$ and $X_1$, then
a Hamiltonian $G$-action on a symplectic DM stack is a presentation
$X_1 \rightrightarrows X_0 \rightarrow \mathcal{X}$ of $\mathcal{X}$
together with a pair of free Hamiltonian $G$-actions on $X_1$ and
$X_0$ (with respect to two possibly degenerate 2-forms) such that the
groupoid structure maps intertwine these  actions. This
is our preferred point of view.
\end{empt}

\begin{proposition}\label{LocallyFreeMoment}
Let $(a : G \times \mathcal{X} \rightarrow \mathcal{X}, 
\ \mu: \mathcal{X} \rightarrow \mathfrak{g}^*)$ 
be a Hamiltonian action of a Lie group $G$ on a symplectic DM stack
$(\mathcal{X}, \omega)$. Then the action is locally free iff the
moment map is regular everywhere (i.e., any value of $\mu$ is
regular).
\end{proposition}
\begin{proof} In a $G$-presentation 
$X_1 \rightrightarrows X_0 \rightarrow \mathcal{X}$
both conditions of the theorem are equivalent 
to the condition that $B_0$ and $A$ are transverse in $TX_0$,
where $A$ is the Lie algebroid of $X_1 \rightrightarrows X_0$
and $B_0$ is as in \ref{InfinitesimalAction}.
\end{proof}
We are now in the position to state and prove the main result of the
section: the DM version of the symplectic
quotient construction.
\begin{theorem}[symplectic reduction]
Let $(a : G \times \mathcal{X} \rightarrow \mathcal{X}, 
\ \mu: \mathcal{X} \rightarrow \mathfrak{g}^*)$ 
be a Hamiltonian action of a compact Lie group $G$ on a symplectic DM
stack $(\mathcal{X}, \omega)$.  Suppose $\xi\in \mathfrak{g}^*$ is a
regular value of $\mu$ which is fixed by the coadjoint action of $G$.
Then 
\[\mathcal{X} /\!\!/\!_{_\xi} G := \mu^{-1} ( \xi )/G
\]
is a DM stack, and $\omega|_{\mu^{-1} ( \xi )}$ descends to a
symplectic form on $\mathcal{X} /\!\!/\!_{_\xi} G$.
\end{theorem}
\begin{proof} By Proposition \ref{RegularValue} $\mu^{-1} ( \xi )$
is a DM stack. Since the coadjoint action of $G$ fixes $\xi$ by
assumption, $G$ acts on $\mu^{-1} ( \xi )$. An argument similar to the
proof of Proposition~\ref{LocallyFreeMoment} shows that the action 
of $G$ on $\mu^{-1} ( \xi )$ is locally free. Hence, by Proposition~\ref{InfinitesimalAction}, $\mu^{-1} ( \xi ) /G$ is a DM stack.

Let $X_1 \rightrightarrows X_0 \rightarrow \mathcal{X}$,
be a $G$-presentation of $\mathcal{X}$. Then 
$\omega = (\omega_1 , \omega_0)$, $\mu = (\mu_1 , \mu_0)$,
and $\mu_1^{-1} ( \xi ) /G \rightrightarrows \mu_0^{-1} ( \xi ) /G 
\rightarrow \mu^{-1} ( \xi ) /G$ is a presentation of $\mu^{-1} ( \xi ) /G$,
where the groupoid maps are induced by those of $X_1 \rightrightarrows X_0$.
Proposition~\ref{InfinitesimalAction} implies that the Lie algebroid $A$ 
of $X_1 \rightrightarrows X_0$ descends to a subbundle 
$\widehat{A}$ of $T(\mu_0^{-1} ( \xi ) /G)$ and $\widehat{A}$
is the Lie algebroid of 
$\mu_1^{-1} ( \xi ) /G \rightrightarrows \mu_0^{-1} ( \xi ) / G$.
Let $B_1 \subset TX_1$ and $B_0 \subset TX_0$ 
be as is in \ref{InfinitesimalAction}.
The moment map equation \eqref{MomentStack} implies that 
$\ker \omega_0 |_{\mu^{-1} ( \xi )} = 
A |_{\mu^{-1} (0)} + B_0 |_{\mu^{-1} (0)}$.
Hence $\omega_0$ descends to a closed 2-form $\widehat{\omega_0}$ 
on $\mu_{0}^{-1} (0) /G$ with the
kernel $\ker \widehat{\omega_0} = \widehat{A}$. Since the groupoid maps
of $\mu_1^{-1} ( \xi ) /G \rightrightarrows \mu_0^{-1} ( \xi ) /G 
\rightarrow \mu^{-1} ( \xi ) /G$ are induced by those of 
$X_1 \rightrightarrows X_0$ we have $t^* \widehat{\omega_0} = s^*
\widehat{\omega_0}$ and so $\widehat{\omega_0}$ defines a symplectic
form on $\mu^{-1} ( \xi ) /G$. 
\end{proof}

\begin{remark}
It is easy to modify the above discussion to describe reduction on a
level $\xi$ which is regular but not coadjoint-invariant. One either
replaces $G$-quotient by the quotient with respect to the stabilizer
of $\xi$ or uses the usual multiplication by the coadjoint orbit
trick.  Nor is the restriction that the group $G$ is compact very
important.  The same result holds for proper actions of non-compact
Lie groups. We concentrate on the simple case to avoid complicating
the notation and to emphasize the stacky features of the reduction.
\end{remark}

\section{Group extension and actions on quotient stacks}\label{sec:3}

A typical example of a DM stack is the quotient $[M/H]$ of a manifold
$M$ by a proper locally free action of a Lie group $H$. This stack is
represented by the Lie groupoid $H \times M \rightrightarrows M$,
where the source map is the projection and the target map is the
action. The corresponding Lie algebroid (as a subbundle of TM) is
spanned by vector fields generating the action of $\mathfrak{h}$, the
Lie algebra of $H$.

\begin{empt}\textbf{Actions on a quotient stack.}\label{MTimesHSection}
\quad 
Suppose we have an exact sequence 
\[
1\to H \to K \to G\to 1
\]
of Lie groups and an action of $K$ on a manifold $M$.  Then the
quotient $G= K/H$ acts on the topological quotient $M/H$.  It is then
reasonable to expect that the extension of $G$ by $H$ also defines an
action of $G$ on the \emph{stack} quotient $[M/H]$.  To define this
action we describe a $G$-atlas of $[M/H]$. Consider a $K$-action
groupoid
\[
K \times (M \times G) \rightrightarrows M \times G
\]
associated to the following $K$-action on $(M \times G)$:
\begin{equation}\label{KMHAction}
k \cdot (m,g) = (k \cdot m, g [k]^{-1}) \ ,
\end{equation}
where $[k]$ is the image of $k$ in $G=K/H$. Note that
$G$ acts on $M\times G$ by left translations on $G$:
$g \cdot (m,g') = (m, gg')$. This $G$-action commutes with the
$K$-action and hence 
$K \times (M\times G) \rightrightarrows M \times G$ is a
$G$-atlas of the stack $[(M \times G) /K]$ with a $G$-action.
Proposition~\ref{ActionExtension} below shows that, in fact,
$[(M \times G) /K]$ is isomorphic to $[M/H]$, where
$H$ acts on $M$ by the restriction of the $K$-action. 
Thus we obtain a $G$-action on $[M/H]$ together with a $G$-atlas
from an action of an extension $K$ on $M$.  
\end{empt}

\begin{proposition}\label{ActionExtension}
The action groupoids
$K \times (M\times G) \rightrightarrows (M \times G)$ and
$H \times M \rightrightarrows M$ are Morita-equivalent,
thus define isomorphic stacks.
\end{proposition}
\begin{proof}
Consider the manifold $M \times K$. 
It has a free $K$-action given by 
\begin{equation}\nonumber
k \cdot (m,k') = (km, k'k^{-1})
\end{equation}
and a commuting free $H$-action given by
\begin{equation}\nonumber
h \cdot (m,k) = (m, hk).
\end{equation}
The map
\begin{equation}\nonumber
\pi_1: \  M \times K \rightarrow M, \quad (m, k) \mapsto  km
\end{equation}
is a principal $K$-bundle and is $H$-equivariant.
The map
\begin{equation}\nonumber
\pi_2: M \times K \rightarrow M \times G, 
\quad (m,k) \mapsto (m, [k])
\end{equation}
is a principal $H$-bundle and is $K$-equivariant, where $K$-action on
$M \times G$ is given by \eqref{KMHAction}.
Thus $ M \xleftarrow{\pi_1} M \times K
\xrightarrow{\pi_2} M\times G$ is a biprincipal bibundle
between the action groupoids $H \times M \rightrightarrows M$ and 
$K \times (M \times G) \rightrightarrows M \times G$.
Thus  the two action groupoids are Morita-equivalent by definition.
\end{proof}

\begin{empt}\textbf{Hamiltonian actions on symplectic quotients.}\quad
We adapt the above constructions to symplectic geometry. Consider a
Hamiltonian action of a compact Lie group $H$ on a symplectic manifold
$(M, \omega_{_M})$ with the moment map $\mu_{_H} : M \rightarrow
\mathfrak{h}^*$.  Let $\xi_{_H}$ be a regular value of
$\mu_{_H}$. Assume $\xi_{_H} \in \mathfrak{h}^*$ is fixed by the
coadjoint action of $H$.  As we remarked above, this is not a very
restrictive assumption since we can always replace $H$ in the
subsequent discussion by the stabilizer of $\xi_{_H}$.  Let $Z =
\mu_{_H}^{-1} (\xi_{_H})$. Then $Z$ is an $H$-invariant manifold and
$[Z / H]$ is a symplectic DM stack, the symplectic quotient $M
/\!\!/_{_{\xi_{_H}}} H:=[Z / H]$.

Next we construct a Hamiltonian action of another Lie group $G$ on the
symplectic stack $[Z/H]$. As in the previous section we consider a Lie
group extension
\[
1 \rightarrow H \rightarrow K \rightarrow G \rightarrow 1 
\]
with the corresponding exact sequences of Lie algebras and dual spaces 
\begin{gather}\nonumber
0 \rightarrow \mathfrak{h} \xrightarrow{i} \mathfrak{k}
\rightarrow \mathfrak{g} \rightarrow 0
\\ \nonumber
0 \leftarrow \mathfrak{h}^* \xleftarrow{i^*} \mathfrak{k}^*
\leftarrow \mathfrak{g}^* \leftarrow 0
\end{gather}
Suppose there is a Hamiltonian action of $K$ on 
$(M, \omega_{_M})$
with the moment map $\mu_{_K}$ such that the restriction
of this action to $H$ is the original $H$-action and the
restriction of $\mu_{_K}$ is $\mu_{_H}$ 
(more precisely $\mu_{_H} = i^* \circ \mu_{_K}$).

The symplectic analogue of the manifold $M \times G$ from
\ref{MTimesHSection} is the symplectic manifold $M \times T^*G$. We
denote by $\lambda^l_{_G}, \lambda^r_{_G}: T^*G \rightarrow
\mathfrak{g}^*$ the projections corresponding to the left and the
right trivializations of $T^*G = \mathfrak{g}^* \times G$
respectively.  The canonical symplectic form on $T^*G$ is given by
$\omega_{_{T^*G}} = d \langle \lambda^l_{_G}, \theta^l_{_G}\rangle$, where
$\theta^l_{_G}$ is the left-invariant Maurer-Cartan form on $G$. With
respect to this symplectic form $-\lambda^l_{_G}$ and $\lambda^r_{_G}$
are the moment maps of the right and the left translations
respectively.

Consider the canonical lift of the action \eqref{KMHAction} of $K$ on 
$M \times G$ to $M \times T^*G$. The lifted action is Hamiltonian 
with respect to the moment map
\begin{equation}
\nu = \mu_{_K} - \lambda^l_{_G} : \ 
M \times T^*G \rightarrow \mathfrak{k}^*
\ ,
\end{equation}
where we think of $\mathfrak{g}^*$ as a subspace of $\mathfrak{k}^*$.

Fix a lift $\xi_{_K} \in \mathfrak{k}^*$ of
$\xi_{_H} \in \mathfrak{h}^*$ and consider the symplectic 
quotient of $M \times T^*G$ at the level 
$\xi_{_K} \in \mathfrak{k}^*$.
We assume $\xi_{_K} \in \mathfrak{k}^*$ is
invariant under the coadjoint action of $K$, 
in particular $Z = \mu_{_K}^{-1} (\xi_{_K})$ is
$K$-invariant.
Again this assumption is not essential but makes statements simpler.
The level set $\nu^{-1} (\xi_{_K})$ can be described as follows
(we use left trivialization of $T^*G$):
\begin{multline}\nonumber
\nu^{-1} (\xi_{_K})  =
\{ (m, \mu_{_K} (m) -  \xi_{_K}, g) \in M \times \mathfrak{g}^* \times G 
= M \times T^*G \ | \ \mu_{_H} (m)  = \xi_{_H} \} \simeq Z \times G
\end{multline}
Hence we get the expected $G$-atlas $Z \times G$ for the DM stack
$[Z/H]$.  Additionally the description of this atlas as 
the moment level set $\nu^{-1} (\xi_{_K})$
provides it with the closed 2-form 
\begin{equation}\nonumber
\omega_{_{Z \times G}} = \omega_{_{M \times T^*G}}|_{\nu^{-1} (\xi_{_K} )} 
= \pi_{_Z}^* \, \omega_{_M} |_{_Z}  + 
d \langle \mu_{_K} \circ \pi_{_Z} - \xi_{_H} , \theta^l_{_G}\rangle ,
\end{equation}
where $\pi_{_Z} : Z \times G \rightarrow Z$ is the projection.
The form $\omega_{_{Z \times G}}$
is $K$-invariant and degenerate precisely along 
$\mathfrak{k}$-action distribution, hence defines a 
symplectic form on the action groupoid 
$K \times (Z\times G) \rightrightarrows Z\times G$.
Moreover the left action of $G$ on $T^*G$ induces a
Hamiltonian action on $Z\times G = \nu^{-1} (\xi_{_K}) \subset
M \times T^*G$ with the Hamiltonian $\eta = \lambda^r_{_G}$
or, explicitly,
\begin{equation}\nonumber
\eta (z, g) = \mathrm{Ad}^*(g) (\mu_{_K} (z) -\xi_{_K}) 
\ : \ Z \times G \rightarrow \mathfrak{g}^* \ .
\end{equation}
Note that different choices of the lift $\xi_{_K}$ of $\xi_{_H}$
correspond to shifts of $\eta$ by a ($G$-invariant) constant.
Also $\eta^{-1} (0) = \mu_{_K}^{-1} (\xi_{_K})$, hence 
the symplectic reduction of the groupoid 
$K \times (Z\times G) \rightrightarrows Z\times G$
with respect to the $G$-action is symplectomorphic to
$K \times \mu_{_K}^{-1} (\xi_{_K})
\rightrightarrows \mu_{_K}^{-1} (\xi_{_K})$,
representing the symplectic quotient of $M$ with respect to $K$.
This is the reduction in stages theorem.

Finally, repeating the proof of Proposition \ref{ActionExtension}
with two commuting free Hamiltonian actions of $K$ and $H$
on $M \times T^*K$ in place of the actions on $M \times K$ one
can easily show that the symplectic action groupoids
$K \times (M\times G) \rightrightarrows (M \times G)$ and
$H \times M \rightrightarrows M$ define isomorphic symplectic stacks.
The modifications required in the proof are similar 
to the above discussion of
the relation between $M \times T^*G$ and $Z \times G$ and
are left to the reader.
\end{empt}
Putting everything together we have the following theorem.

\begin{theorem}\label{ActionSymplecticQuotientTheorem}
Let $1 \rightarrow H \rightarrow K \rightarrow G
\rightarrow 1$ be a sequence of compact Lie groups and 
$i^*:\mathfrak{k}^*\to\mathfrak{h}^*$ the corresponding canonical
projection.  Suppose there is a Hamiltonian action of $K$ on a
symplectic manifold $(M, \omega_{_M} )$ with a moment map $\mu_{_K}: M
\to \mathfrak{k}^*$.  Suppose $\xi_{_K} \in
\mathfrak{k}^*$ is $K$-invariant and suppose 
$\xi_{_H}:= i^* (\xi_{_K})$ is a regular value of the $H$-moment map
$\mu_{_H}:= i^*\circ \mu_{_K}$.  Then
\begin{enumerate}
\item
$[\mu_{_H}^{-1} (\xi_{_H})/H]$ is a symplectic DM stack with the
symplectic form induced by the restriction of $\omega_{_M}|_{\mu_H^{-1}
(\xi_{_H})}$;
\item\label{HHamiltonianActionTheorem} 
There is a Hamiltonian action of $G$ on $[\mu_{_H}^{-1} (\xi_{_H})/H]$ 
induced by the action of $K$ on $M$;
\item 
Assume $\xi_{_K}$ is a regular value of $\mu_{_K}$. Then the
symplectic reduction 
$[\mu_{_{_H}}^{-1} (\xi_{_H})/H] /\!\!/_{_0} G$
at $0$ of $[\mu_{_H}^{-1} (\xi_{_H})/H]$ with respect to the $G$-action
defined in (\ref{HHamiltonianActionTheorem}) is a 
symplectic DM stack isomorphic to $[\mu_{_K}^{-1} (\xi_{_K})/K]$:
\[
[\mu_{_H}^{-1} (\xi_{_H})/H] /\!\!/_{_0} G \simeq [\mu_{_K}^{-1}
(\xi_{_K})/K].
\]
\end{enumerate}
\end{theorem}

\section{Symplectic toric DM stacks}\label{sec:4}

\begin{definition} \emph{A symplectic $G$-toric DM stack} is 
a symplectic DM stack $\mathcal{X}$
with a Hamiltonian action of a compact torus $G$, such that
\begin{itemize}  
\item $\dim \mathcal{X} = 2 \dim G$ and 
\item the action of $G$ on the coarse moduli space of 
$\mathcal{X}$ is effective.
\end{itemize}
\end{definition}
This is a natural definition of a toric object in the context of
symplectic DM stacks.

\begin{empt}\textbf{Finite extension of a torus.}
As an example of an application of 
Theorem~\ref{ActionSymplecticQuotientTheorem} we construct a symplectic toric
DM stack as symplectic quotients of $\mathbb{C}^N$ (a generalization
of Delzant's construction \cite{Delzant1988} of symplectic toric
manifolds).  We start with an extension
\[
1\to \Gamma \to \widehat{T^N} \to T^N \to 1
\]
of the $N$-dimensional compact torus $T^N :=
\mathbb{R}^N/\mathbb{Z}^N$ by a finite group $\Gamma$. Then for any
closed subgroup $A$ of $T^N$ we have an extension $\widehat{A}$ by
$\Gamma$.  The standard action of $T^N$ on the symplectic space
$(\mathbb{C}^N, \omega_{_{\mathbb{C}^N}} = \sqrt{-1} \sum_j dz_j
\wedge d \bar{z}_j)$ gives rise to an
(ineffective) Hamiltonian action of $\widehat{T^N}$ on $(\mathbb{C}^N, \omega_{_{\mathbb{C}^N}})$ 
with the ``same'' moment map 
$\mu: \mathbb{C}^N \rightarrow (\mathbb{R}^N)^*$. 
We would like to apply Theorem \ref{ActionSymplecticQuotientTheorem}
which requires a suitable ($\widehat{T^N}$-invariant)
choice of the moment map level $\xi$. The following Lemma ensures that any level of the moment map works fine in our situation.
\end{empt}
\begin{lemma}
Let  $1\to \Gamma \to \widehat{T} \to T \to 1$ an extension of  a connected abelian Lie group $T$ by a finite group $\Gamma$.  Then the coadjoint
representation of $\widehat{T}$ is trivial.
\end{lemma}
\begin{proof}
Denote the connected component of the identity of 
$\widehat{T}$ by $\widehat{T}_0$.  
Since $\widehat{T}_0$ is a connected cover of an abelian
Lie group, it is an abelian Lie group and its adjoint action is
trivial. Hence it is enough to show that the adjoint action of
$\Gamma$ is trivial. Since $\Gamma$ is finite, for any $\gamma\in \Gamma$ $X \in \widehat{\mathfrak{t}}$, and $t\in \mathbb{R}$, we have
\[
\exp (tX) \gamma \exp (-tX) = \exp (0X) \gamma \exp (-0X) = \gamma.
\]
Hence 
\[
\gamma\exp (tX) \gamma ^{-1} = \exp (tX)
\]
for all $t\in\R$.  Taking derivatives of both sides with respect to
$t$ at $t=0$ we get $Ad(\gamma) X = X$. $\Box$
\end{proof}

\begin{thm}\label{thm:5.4}
Let $1\to \Gamma \to \widehat{T^N} \to T^N \to 1$ be an extension of
the standard $N$-dimensional torus $T^N= \R^N/\Z^N$ by a finite group
$\Gamma$, $A<T^N$ a closed subgroup and $\widehat{A}<\widehat{T^N}$
the corresponding subgroup of $\widehat{T^N}$.  Let $a\in
\mathfrak{a}^* $ be a regular value of the $A$-moment map
$\mu_{\widehat{A}}:\C^N\to \mathfrak{a}^*$.  Then the stack quotient
\[
\C^N/\!\!/_a \hat{A} = [\mu_{\widehat{A}}^{-1} (a)/\widehat{A}]
\]
is a symplectic toric $G$-manifold, where $G=
\widehat{T^N}/\widehat{A} = T^N/A$.
\end{thm}
\begin{remark}
The same stack can be obtained as the reduction of the symplectic DM
stack $\mathbb{C}^N \! \times \! B \Gamma$ with respect to the
diagonal Hamiltonian action of $A$, where the action of $A$ on
$B\Gamma$ is defined by the extension $1 \rightarrow \Gamma
\rightarrow \widehat{A} \rightarrow A \rightarrow 1$:
\[
\mathbb{C}^N/\!\!/\!_{_a} \widehat{A} \cong 
(\mathbb{C}^N \! \times \! B \Gamma ) /\!\!/\!_{_a} A \ .
\]
This is the point of view taken in the Introduction.
\end{remark}
\begin{proof}[Proof of Theorem~\ref{thm:5.4}. ]
By Theorem~\ref{ActionSymplecticQuotientTheorem} with
$K=\widehat{T^N}$, $H = \widehat{A}$, $G = \widehat{T^N}/\widehat{A} =
T^N/A$ and $(M, \omega) = (\mathbb{C}^N,\omega_{_{\mathbb{C}^N}})$ the
symplectic quotient
\[
\mathbb{C}^N/\!\!/\!_{_a} \widehat{A} := 
[\mu_{\widehat{A}}^{-1} (a)/\widehat{A}]
\]
is a symplectic DM stack with a Hamiltonian action of the
torus $G$ for any 
regular value $a \in \mathfrak{a}^*$ of the moment map 
$\mu_{\widehat{A}}$. 
We have $\dim \mathbb{C}^N/\!\!/\!_a \widehat{A} = 2N - \dim A - \dim
A = 2 \dim G$.  So to prove that the symplectic DM stack
$\mathbb{C}^N/\!\!/\!_{_a}
\widehat{A}$ is $G$-toric it remains to check that the action of $G$
on its coarse moduli space $ \mu_{\widehat{A}}^{-1} (a)/\widehat{A}$ is
effective.

We first argue that $Z \cap (\mathbb{C}^\times)^N \neq \emptyset$,
where $Z=\mu_{_H}^{-1} (a)$. The moment map $\mu_{_K} : \mathbb{C}^N
\rightarrow (\mathbb{R}^N)^* = \mathfrak{k}^*$ for the action of
$K=\widehat{T^N}$ on $\mathbb{C}^N$ is given by $\mu_{_K} (z_1,
\ldots, z_N) =
\sum |z_j|^2 e_j^*$, where $\{ e_j^* \}$ is the basis of 
$(\mathbb{R}^N)^*$ dual to the standard basis $\{e_j\}$ 
of $\mathbb{R}^N = \mathfrak{k}$ 
(which is also a basis of $\mathbb{Z}^N$).  
The image $\mu_{_K} (\mathbb{C}^N)$ is the orthant
\begin{equation}\nonumber
(\mathbb{R}^N)^*_+ := 
\{ \eta \in (\mathbb{R}^N)^* \mid \langle \eta, e_j \rangle 
\geq 0, 1 \leq j \leq N \}
\end{equation}
with $\mu_{_K}$ mapping $(\mathbb{C}^\times)^N$ 
to the interior of the orthant.  Since
$\mu_{_{_H}}$ is the restriction of $\mu_{_K}$, we have
\begin{equation}
Z \cap (\mathbb{C}^\times)^N \neq \emptyset 
\quad \Leftrightarrow \quad
V_a \cap \textrm{interior} ((\mathbb{R}^N)^*_+) \neq \emptyset
\ ,
\end{equation}
where the affine subspace 
$V_a \subset (\mathbb{R}^N)^* =\mathfrak{k}^*$ is 
the preimage of $a \in \mathfrak{h}^*$. If $\tilde{a} \in V_a$
then $V_a = \tilde{a} + \mathfrak{h}^\perp$, 
where $\mathfrak{h}^\perp$ denotes the annihilator
of $\mathfrak{h}$ in $(\mathbb{R}^N)^*$.

The faces of the orthant $(\mathbb{R}^N)^*_+$ 
are images under $\mu_{_K}$ of coordinate subspaces of the form
\begin{equation}\nonumber
\mathcal{S} =\mathcal{S} (i_1, \ldots, i_n) :=
\{ z_{i_1} = 0, \ldots, z_{i_n} = 0 \} 
\end{equation}
for some subset $\{ i_1, \ldots, i_n \} \subset \{ 1,\ldots, N \}$.
The subspace $\mathcal{S}$ above is precisely 
the fixed point set of the subtorus
\begin{equation}\nonumber
S = \bigl\{ ( \lambda_1, \ldots, \lambda_N) \in T^N \mid 
\lambda_j = 1 \textrm{ for } i\not\in \{ i_1, \ldots, i_n \} \bigr\}
\end{equation}
with Lie algebra $\mathfrak{s}$. 
Since the action of $H$ on $Z$ is locally
free, $Z \cap \mathcal{S} \neq \emptyset$ 
implies that $\mathfrak{s} \cap \mathfrak{h} = 0$.
Hence $\mathfrak{s}^\perp + \mathfrak{h}^\perp = (\mathbb{R}^N)^*$.
We conclude that the affine
plane $V_a$ intersects the faces of the orthant 
$(\mathbb{R}^N)^*_+$
transversely.  Therefore $V_a$ contains points in the interior of $(\mathbb{R}^N)^*_+$, hence $Z \cap (\mathbb{C}^\times)^N
\not = \emptyset$, and there is a point $z \in Z$ on which 
$T^N$ acts freely.

An element $g\in G$ acts trivially on the coarse moduli space of 
the stack $[ Z / H] = [(Z \times G)/\widehat{T^n}]$ if for any
$(z, g') \in Z \times G$ there is 
$\widehat{x} \in \widehat{T^N}$ such that
\begin{equation}\nonumber
\widehat{x} \cdot (z,g') = g \cdot (z, g') \ , 
\end{equation}
that is,
\begin{equation}\label{StabEquation}
(\widehat{x} z, g' [\widehat{x}]^{-1} ) = (z, g g') \ ,
\end{equation}
where $[\widehat{x}] \in G = \widehat{T^n}/\widehat{A}$ 
is the class of $\widehat{x} \in \widehat{T^n}$.
Now take $g'=e_{_G}$, the identity of $G$, and $z \in Z$ an
element on which $T^N$ acts freely. Then \eqref{StabEquation}
implies
\begin{equation}\nonumber
\widehat{x} \in \Gamma \ , \quad g=[\widehat{x}]=e_{_G} 
\end{equation}
and, since $\Gamma \subset \widehat{A}$, we have $g=e_{_G}$.
Hence the action of $G$ on the coarse moduli space $[Z / \widehat{A}]$ 
is effective, and $\mathbb{C}^N/\!\!/\!_{_a} \widehat{A}$
is a symplectic toric DM stack.
\end{proof}

\end{document}